\documentclass{amsproc}

\newtheorem{theorem}{Theorem}[section]
\newtheorem{lemma}[theorem]{Lemma}

\theoremstyle{definition}
\newtheorem{definition}[theorem]{Definition}

\theoremstyle{remark}

\numberwithin{equation}{section}

\begin{document}

\title{The Parameter Rigid Flows on Orientable 3-Manifolds}

\author{Shigenori Matsumoto}
\address{Department of Mathematics, College of
Science and Technology, Nihon University, 1-8-14 Kanda, Surugadai,
Chiyoda-ku, Tokyo, 101-8308 Japan
}
\curraddr{Department of Mathematics, College of
Science and Technology, Nihon University, 1-8-14 Kanda, Surugadai,
Chiyoda-ku, Tokyo, 101-8308 Japan}
\email{matsumo@math.cst.nihon-u.ac.jp
}
\thanks{The author was partially supported by Grant-in-Aid for
Scientific Research (A) No.\ 17204007.}

\date{January 10, 2009 }


\begin{abstract}
A flow defined by a nonsingular smooth vector field $X$ 
on a closed manifold $M$ is said to be
parameter rigid if given any real valued smooth function $f$ on $M$, 
there are a smooth funcion $g$ and a constant $c$ such that
 $f=X(g)+c$ holds.
We show that the parameter rigid flows on closed orientable 3-manifolds are 
smoothly conjugate to Kronecker flows on the 3-torus with badly approximable 
slope.

\end{abstract}

\maketitle

\section{Introduction}

Throughout this paper we work in the $C^\infty$-category:
any manifold, function, diffeomorphism, form, vector field e.t.c.\
are to be of class $C^\infty$.
Let $X$ be a nonsingular vector field on a closed manifold $M$
that defines a flow $\varphi^t$.

\begin{definition} \label{def1}
\rm  The flow $\varphi^t$  
is called {\em parameter rigid} 
if for any function $f$ on $M$, there are a function $g$ and a constant
$c$ such that $f=X(g)+c$ holds.
\end{definition}

It is well known and easy to show
that the parameter rigidity is equivalent to the following
property: if $\psi^t$ be another nonsingular flow 
defined 
by a vector field $fX$, where $f$ is a nowhere vanishing function,
then there are an orbit preserving diffeomorphism $F$ of $M$ and a nonzero
constant $c$ such that
$$
\psi^t(F(x))=F(\varphi^{ct}(x)).
$$

The only known examples of parameter rigid flows are Kronecker flows on tori with
badly approximable (sometimes called Diophantine or non-Liouville) slope,  and A. Katok has conjectured that in fact they are the all
(\cite{K}).
In this paper we show a partial result supporting this conjecture.

\begin{theorem} \label{t1}
A parameter rigid flow on a closed orientable 3-manifold  
is smoothly conjugate to
a linear flow on the 3-torus with badly approximable slope.
\end{theorem}

At Paulfest, A. Kocsard has announced the same result (\cite{Ko}).

The method of this paper cannot be applied to the nonorientable
3-manifolds. The difficulty lies in showing that the lift of
a parameter rigid flow to the orientable double cover is
again parameter rigid.

Thanks are due to the unanimous referee, whose valuable comments
are helpful for the shorter and clearer arguments. 

\section{General properties of parameter rigid flows}

Here we collect some basic facts needed in
 the proof of Theorem \ref{t1}.
Let $\varphi^t$ be a parameter rigid flow on a closed
$(n+1)$-dimensional manifold, defined by a nonsingular
vector field $X$. 
\\[2mm]
(1) {\em The flow $\varphi^t$ is uniquely ergodic, leaves a 
volume form $\Omega$ invariant, and hence is minimal.}

\smallskip
Indeed the Birkhoff average of any smooth function tends to
a constant, which is enough for the unique ergodicity, since
the smooth functions are dense in the space of continuous
functions. For the second statement, let $\Omega_0$ be
an arbitrary volume form and define a function $f$
by ${\mathcal L}_X\Omega_0=f\Omega_0$. Then $\Omega=e^{-g}\Omega_0$
is the desired form, where $g$ is the function obtained
by Definition \ref{def1}.
\\[2mm]
(2) {\em The function $g$ in Definition \ref{def1} is unique up
to a constant sum, and the constant $c$ is given by}
$c=\int_Mf\Omega$.

\smallskip
In fact if $X(h)$ is constant, then it should be 0,
and the minimality of the flow implies that $h$ is constant.
\\[2mm]
(3) {\em The vector space  $\Lambda^n(X)$ consisting
of $n$-forms $\omega$ such that
$i_X\omega=i_Xd\omega=0$ is one dimensional, spanned by
$i_X\Omega$.}

\smallskip
Indeed $i_X\Omega$ belongs to $\Lambda^n(X)$ and any $n$-form
in $\Lambda^n(X)$ is a function multiple of $i_X\Omega$.
Taking the Lie derivative, one can show the function is
constant.

\smallskip
A 1-form $\alpha$ is called {\em normal} if $\alpha(X)$ is constant.
The {\em normalization} $\alpha$ of any 1-form $\alpha'$ is defined
to be $\alpha=\alpha'-dg$, where $g$ is a function (unique
up to constant sum) such that
$\alpha'(X)=X(g)+c$.
\\[2mm]
(4) {\em A closed normal 1-form $\alpha$ is invariant
by the flow
i.\ e.\ ${\mathcal L}_X\alpha=0$. 
By the minimality of the flow, it is either identically zero or nonsingular.}

\smallskip
Let $\Lambda^1(X)$ be the space of closed normal 1-forms and let
$\epsilon:\Lambda^1(X)\to H^1(M;{\mathbb R})$ be the map assigning the
cohomology class to each form.
\\[2mm]
(5){\em The homomorphism $\epsilon$ is an isomorphism.}

\smallskip
Indeed the normalization of each closed form belongs to $\Lambda^1(X)$,
showing the surjectivity of $\epsilon$. On the other hand an 
exact normal form
$dg$ is identically zero, since  if $X(g)$ is a constant, then 
$g$ is a constant.

\section{Proof of the main theorem}

Let $\varphi^t$ be a parameter rigid flow
defined by a vector field $X$ on a closed orientable 
3-manifold $M$. We shall prove Theorem \ref{t1}
dividing into cases.

\bigskip

\noindent
{\bf Case 1}. $H^1(M;{\mathbb R})\neq 0$.

\smallskip

Let $\alpha\in\Lambda^1(X)$ be a closed normal 1-form representing
an integral class. Then the equation $\alpha=0$ defines  a fibration of $M$
over the circle. The constant $\alpha(X)$ cannot be 0, since
the flow is minimal. Thus the flow has  
a global cross section,
say $\Sigma$. The first return map of $\Sigma$ must be minimal, 
and especially it does not admit any periodic point.
Then by a theorem of Jiang \cite{J}, one can show
that  $\Sigma$ is  diffeomorphic to the 2-torus. 
Now the first return map is cohomologically rigid
in the sense of \cite{LS} and is shown in that paper
to be conjugate to a translation by a badly approximable vector.
We have done with this case.

\bigskip

\noindent
{\bf Case 2}. $H^1(M;{\mathbb R})=0$.

\smallskip

Let $\Omega$ be the volume form which is left invariant by $X$.
Then since
$$
{\mathcal L}_X\Omega=d\iota_X\Omega=0,
$$
there is a 1-form $u$ such that
$
\iota_X\Omega=du
$. Taking the normalization of the previous section, one may assume
that $u$ is normal, i.\ e.\
$u(X)=c_1$ is a constant. Then since $\iota_X(u\wedge du)=c_1du=\iota_X(c_1\Omega)$,
we have $u\wedge du=c_1\Omega$.

\bigskip

\noindent
{\bf Case 2.1}. $c_1\neq 0$.

\smallskip

In this case the vector field $c_1^{-1}X$ generates the Reeb flow of a contact form $u$. The solution of the Weinstein conjecture
in \cite{T} shows that the flow admits a closed orbit, contrary to
the minimality.

\bigskip

\noindent
{\bf Case 2.2}. $c_1=0$.

\smallskip
We have $\iota_X\Omega=du$ and $u(X)=0$.
First of all notice that $u$ is nonsingular. Indeed we have 
$\mathcal L_Xu=0$, that is, $u$ is invariant by the flow $\varphi^t$.
By the minimality of the flow $\varphi^t$, vanishing of $u$ at some point
would imply that $u$ is identically zero, which is not the case since
$du$ is nonsingular. As noted before, we have $u\wedge du=0$,
that is, the 1-form $u$ is integrable,
and $du=\eta'\wedge u$ for some 1-form $\eta'$. 
Notice that $\eta'(X)=0$.
We get
$$
0=d(\iota_Xdu)=\mathcal L_Xdu=\mathcal L_X\eta'\wedge u+\eta'\wedge 
\mathcal L_Xu.
$$
Since $\mathcal L_Xu=0$, we have $\mathcal L_X\eta'\wedge u=0$. That is,
$\mathcal L_X\eta'=f_2u$ for some function $f_2$. Write $f_2
=X(g_2)+c_2$ and let $\eta=\eta'-g_2u$. Then we have
$$
du=\eta\wedge u , \ \ \ \eta(X)=0 \ \ \ {\rm and} \ \ \ \mathcal L_X\eta=c_2u.
$$

\bigskip

\noindent
{\bf Case 2.2.1}. $c_2=0$.

\smallskip

In this case we have $d\eta\in \Lambda^2(X)$, and thus by (3) of
the previous section, $d\eta=r\iota_X\Omega$ for some constant $r$.
Since $\eta\wedge u=du=\iota_X\Omega$ is nonvanishing, $ru-\eta$ is a nonzero element of $\Lambda^1(X)\cong H^1(M;{\mathbb R})$,
contrary to the assumption of Case 2.

\bigskip

\noindent
{\bf Case 2.2.2}. $c_2\neq 0$.

\smallskip
Changing $X$ and $u$ by a scalar multiple at the same time one may assume that $\mathcal L_X\eta= -2u$ and 
still $du=\iota_X\Omega$.
In summary, there are two 1-forms $u$ and $\eta$ such that
$$
du=\eta\wedge u=\iota_X\Omega, \ \ u(X)=\eta(X)=0, \ \ \iota_Xd\eta=-2u.
$$ 

Since $\eta\wedge u$ is nonvanishing, that is, $\eta$ and
$u$ are linearly independent everywhere, there is a 1-form
$\sigma'$ such that $\Omega=\eta\wedge u\wedge\sigma'$.
Then the triplet $\langle \eta,u,\sigma'\rangle$ is a
basis of the space of 1-forms as a module over the ring of
functions, and
likewise $\langle\eta\wedge u,u\wedge\sigma',\sigma'\wedge\eta\rangle
$ is a basis of the space of 2-forms.

Note that $\sigma'(X)=1$
since $\iota_X\Omega=du$.
Now we have:
$$
0=\mathcal L_X\Omega=\mathcal L_X(du\wedge\sigma')=\mathcal L_Xdu\wedge \sigma'+du\wedge\mathcal L_X\sigma'.
$$
But $\mathcal L_Xdu=0$, and thus we have $du\wedge\mathcal L_X\sigma'=0$.
Thus one can write
$$
\mathcal L_X\sigma'=f_3\eta+f_4u.
$$
Then there are functions $g_3$ and $g_4$ such that
$$
f_3=X(g_3)+c_3, \ \ \ f_4+2g_3=X(g_4).
$$
(In the last expression, we do not need a constant, since we can alter $g_3$
by a constant summand.)
Now computation shows that for $\sigma=\sigma'-g_3\eta-g_4u$, we have
$\mathcal L_X\sigma=c_3\eta$. Summing up, we have obtained
$$
\Omega=\eta\wedge u\wedge\sigma, \ \ \ \sigma(X)=1, \ \ \
\mathcal L_X\sigma=c_3\eta.
$$

We prepare a useful lemma.

\begin{lemma} \label{l1}
If $\mathcal L_Xw=a\Omega$ for some 3-form $w$ and a constant $a$,
then $a=0$ and the form $w$ is invariant by $X$.
\end{lemma}

\noindent
\begin{proof}
The proof is immediate by taking the integral over $M$.
\end{proof}

\bigskip

We are going to show that in fact the manifold $M$ is
a quotient of a 3-dimensional Lie group. For this we need to compute
$d\eta$ and $d\sigma$. First of all let
$$
d\eta=f_5\eta\wedge u+f_6\eta\wedge\sigma+f_7u\wedge\sigma.
$$
Then since $\iota_Xd\eta=-f_6\eta-f_7u=-2u$, we have $f_6=0$ and $f_7=2$,
that is, $d\eta=f_5\eta\wedge u +2u\wedge \sigma$.
Now
\begin{align*}
\mathcal L_X(\sigma\wedge d\eta)
&=c_3\eta\wedge d\eta+\sigma\wedge \mathcal L_Xd\eta\\
&=c_3\eta\wedge d\eta+\sigma\wedge d(\iota_Xd\eta)\\
&=c_3\eta\wedge d\eta+\sigma\wedge (-2du)\\
&=2c_3\eta\wedge u\wedge \sigma -2\sigma\wedge \eta\wedge u\\
&=2(c_3-1)\eta\wedge u\wedge\sigma.
\end{align*}

By lemma \ref{l1}, we have
$$
c_3=1 \ \ \ {\rm and} \ \ \ \mathcal L_X(\sigma\wedge d\eta)=0.
$$

On the other hand, we have
$$
\mathcal L_X(\sigma\wedge d\eta)=\mathcal L_X(f_5\sigma\wedge \eta\wedge u)=\mathcal L_X(f_5\Omega)
=X(f_5)\Omega.
$$
Thus  by Lemma \ref{l1}, we have $X(f_5)=0$, that is,  $f_5$ is a constant, say $c_5$.

In summary we obtained:
$$
\mathcal L_X\sigma=\eta, \ \  \ d\eta=c_5 \eta\wedge u +2 u\wedge \sigma.
$$

\bigskip

An unknown constant $c_5$ will be shown to be zero in the way of
computing $d\sigma$. Let
$$
d\sigma=f_8\eta\wedge u+ f_9\eta\wedge\sigma +f_{10}
u\wedge\sigma.
$$

Since $\eta=\iota_Xd\sigma=-f_9\eta-f_{10}u$, we have $f_9=-1$ and $f_{10}=0$.
That is, $d\sigma=f_8\eta\wedge u- \eta\wedge\sigma$.
Then we have
\begin{align*}
\mathcal L_X(\sigma\wedge d\sigma) 
&=\eta\wedge d\sigma +\sigma\wedge \mathcal L_Xd\sigma=
0+\sigma\wedge d\eta\\
&=\sigma\wedge(c_5\eta\wedge u+2u\wedge\sigma)
=c_5\sigma\wedge \eta\wedge u=c_5\Omega.
\end{align*}
Again by Lemma \ref{l1} we conclude that $c_5=0$. 

On the other hand,
$$
\mathcal L_X(\sigma\wedge d\sigma)=\mathcal L_X(f_8\sigma\wedge \eta\wedge u)=\mathcal L_X(f_8\Omega)
=X(f_8)\Omega.
$$
This implies $X(f_8)=0$.  That is, $f_8$ is a constant $c_8$.

Summing up, one gets
$$
du=\eta\wedge u, \ \ d\eta=2u\wedge\sigma, \ \ d\sigma=c_8\eta\wedge u -\eta
\wedge \sigma.
$$

Letting $\hat\sigma=\sigma-c_8u$, we obtain a final conclusion.

\begin{lemma} \label{l2}
On the manifold $M$, there are three 1-forms $\eta$, $u$ and
$\hat\sigma$ such that
\begin{align*}
& \Omega=\eta\wedge u\wedge\hat\sigma, \\
& d\eta=2u\wedge\hat\sigma, \ \ \ du=\eta\wedge u, \ \ \
d\hat\sigma=-\eta\wedge\hat\sigma, \\
& \eta(X)=0, \ \ \ u(X)=0, \ \ \ \hat\sigma(X)=1.
\end{align*}
\end{lemma}

This lemma says that the manifold $M$ is the quotient space
of the universal cover of the Lie group ${\rm SL}(2,{\mathbb R})$
by a cocompact lattice, and the vector field $X$ generates the horocycle flow.
But the horocycle flow is shown not to be 
parameter rigid by \cite{FF}. So this case also leads to a contradiction, and we have done with the proof of Theorem \ref{t1}.

\end{document}